\documentclass[12pt]{article}
\usepackage{amsmath, amsthm, amscd, amsfonts, amssymb, graphicx, color}
\usepackage{epstopdf}
\usepackage{subfigure}
\usepackage{setspace}
\usepackage[left=1in, right=1in, bottom=1.4in]{geometry}
\usepackage{bm}
\usepackage{cite}
\newtheorem{theorem}{Theorem}[section]

\newtheorem{lemma}[theorem]{Lemma}

\newtheorem{definition}{Definition}

\title{ \bf  Some optimal inequalities for $\alpha$-harmonic functions estimated by their boundary functions$^{\ast}$ }
\date{}
 \begin{document}

\maketitle
\renewcommand{\thefootnote}{\fnsymbol{footnote}}
\noindent{\footnotesize\rm{ \begin{center}\bf Bo-Yong Long  \end{center}\
 \begin{center}  School of Mathematical Sciences, Anhui University, Hefei  230601, China
\end{center}}
\footnotetext{\hspace*{-5mm}
\begin{tabular}{@{}r@{}p{16.0cm}@{}}
&$^{*}$Supported by NSFC (No.12271001), Natural Science Foundation of Anhui Province (2308085MA03), and Excellent University Research and Innovation Team in Anhui Province (2024AH010002), China.
\\
&E-mail:  boyonglong@163.com

\end{tabular}}

\begin{quote}
{\small \noindent {\bf Abstract}:\   The solutions of a kind of second-order homogeneous  partial differential  equation  are called (real kernel) $\alpha$-harmonic functions.   The $\alpha$-harmonic functions and their first-order partial derivative functions on unit disk are  estimated using the $L^{p}$ norm of the boundary functions of the $\alpha$-harmonic functions. A series of inequalities are obtained. In addition, when the $\alpha$-harmonic functions are quasiconformal,   their first-order partial derivative functions are estimated by the arc length of the domain boundary and the Lipschitz constant of the boundary functions. All of the inequalities obtained in this article are optimal or asymptotically optimal.
}\\
  {{\small \noindent{\bf Keywords}: $\alpha$-harmonic functions;  boundary values;  pointwise estimate; quasiconformal mappings }
 }
\end{quote}

\section{Introduction and main results}
\setcounter{equation}{0}
\hspace{2mm}

{\bf 1.0  $\alpha-$harmonic functions}

Let $\mathbb{D}=\{z:|z|<1\}$ be the  open unit disk  and $\mathbb{T}=\{z:|z|=1\}$ unit circle.
For $\alpha\in\mathbb{R}$ and  $z\in \mathbb{D}$, let

$$T_{\alpha}=-\frac{\alpha^{2}}{4}(1-|z|^{2})^{-\alpha-1}+\frac{\alpha}{2}(1-|z|^{2})^{-\alpha-1}\left(z\frac{\partial}{\partial z}
+\bar{z}\frac{\partial}{\partial\bar{z}}\right)+\frac{1}{4}(1-|z|^{2})^{-\alpha}\triangle$$
be a second-order elliptic  partial differential operator, where $\triangle$ is the usual complex Laplacian operator
$$\triangle:=4\frac{\partial^{2}}{\partial z \partial \bar{z}}=\frac{\partial^{2}}{\partial x^{2}}+\frac{\partial^{2}}{\partial y^{2}}. $$
The corresponding homogeneous differential equation is
\begin{align}\label{1.1}
T_{\alpha}(u)=0 \quad \mbox {in}  \,\,\mathbb{D}
\end{align}
and the associated the Dirichlet boundary value problem is as follows
\begin{equation} \label{1.2}\;  \left\{
\begin{array}{rr}
T_{\alpha}(u)=0 &\quad \mbox {in}  \,\,\mathbb{D}, \\
u=f &\quad \mbox {on}  \,\,\mathbb{T}.
\end{array}\right.
\end{equation}
Here, the boundary data $f\in\mathfrak{D}'(\mathbb{T})$
is a distribution on the boundary $\mathbb{T}$ of $\mathbb{D}$, and the boundary condition in (\ref{1.2}) is interpreted in the distributional sense that $u_{r}\rightarrow f$ in $\mathfrak{D}'(\mathbb{T})$
as $r\rightarrow1^{-}$, where
\begin{equation*}
u_{r}(e^{i\theta})=u(re^{i\theta}), \quad e^{i\theta}\in\mathbb{T},
\end{equation*}
for $r\in[0,1)$. By the results of  \cite{Olofsson2014}, we know that if $\alpha>-1$ and $u\in\mathcal{C}^{2}(\mathbb{D})$ satisfy (\ref{1.2}), then  the function $u$ has the form of Poisson type integral
\begin{equation}\label{1.3}
u(z):=P_{\alpha}[f](z)=\frac{1}{2\pi}\int_{0}^{2\pi}K_{\alpha}(ze^{-it)})f(e^{it})dt, \quad \mbox{for}\,\, z\in\mathbb{D},
\end{equation}
where
 \begin{equation}\label{1.4}
 K_{\alpha}(z)=c_{\alpha}\frac{(1-|z|^{2})^{\alpha+1}}{|1-z|^{\alpha+2}},
 \end{equation}
$c_{\alpha}=\Gamma^{2}(\alpha/2+1)/\Gamma(1+\alpha)$ and $\Gamma(\cdot)$  is the  Gamma function.

It was shown in \cite{Olofsson2014} that if $\alpha\leq -1$, $u\in\mathcal{C}^{2}(\mathbb{D})$ satisfy (\ref{1.1}) and  the boundary limit  $f=\lim_{r\rightarrow 1^{-}}u_{r}$ exists in $\mathfrak{D}'(\mathbb{T})$, then
  $u(z)=0$ for all $z\in\mathbb{D}$.  That is to say, the case of $\alpha\leq -1$ is trivial. So, in the following we will only focus on the case of $\alpha>-1$.
We call the solutions of (\ref{1.1})  with $\alpha>-1$
 real kernel $\alpha$-harmonic functions \cite{Long2022Proc, Long2024, Long2021Filomat}, or  $\alpha$-harmonic   functions.

On the one hand, if we take $\alpha=2(p-1)$, then $u$ is polyharmonic (or $p-$harmonic), where $p\in\{1,2,...\}$ (cf.\cite{Du2012, Chens2021, Lip2022, Lium2024, Iwaniec2019}).
In particular, if $\alpha=0$, then $u$ is harmonic (cf.\cite{Kalaj2015, Kalaj2019, Bshouty2018}).  Thus, $u$ is a kind of generalization of classical harmonic mappings. On the other hand, $\alpha$-harmonic functions themselves have also been generalized to so-called $(\alpha, \beta)$-harmonic functions or $(p, q)$-harmonic functions \cite{Klintborg2021, Arsenovic2024}.

{\bf 1.1  Estimating $\alpha$-harmonic functions using boundary functions}

 Let $u$ be $\alpha$-harmonic function defined on the unit disk $\mathbb{D}$ and $f$ the corresponding boundary function. In \cite{Kalaj2024}, the author obtains some sharp or asymptotically  sharp estimates of the type $|u(z)|\leq g(r) \|f\|_{L^{p}(\mathbb{T)}}$ and $|Du(z)|\leq h(r)\|f\|_{L^{p}(\mathbb{T)}}$. There are also similar  results of $\alpha$-harmonic functions \cite{Chenj2023}. For similar results but for the harmonic mappings in the unit disk see  \cite{Zhu2021, Colonna1989}. For the multidimensional setting see \cite{Khavinson1992, Kresin2010, Liucw2021, Melentijevic2019Adv}. For hyperbolic harmonic or so-called $n$-harmonic mappings in the unit ball see \cite{Chenj2023PA, KhalfallahA2022}.

 In this section, we will conduct further research on such issues. Let $M_{p}(r, g)$ denote integral means of $g$. Firstly, we arrive at the following concise and interesting conclusion.

\begin{theorem}\label{Th1.6}
Let $u(z)=P_{\alpha}[f](z)$ be $\alpha$-harmonic functions with $f\in L^{p}(\mathbb{T})$, $p\geq 1$ and $z=re^{i\theta}\in\mathbb{D}$. Then the following sharp inequality
\begin{align}\label{1.11}
M_{p}(r,u)\leq c_{\alpha}F(-\alpha/2, -\alpha/2; 1; r^{2})\|f\|_{L^{p}(\mathbb{T)}}
\end{align}holds. In particular, it holds that
\begin{align}\label{1.12}
\|u\|_{H^{p}(\mathbb{T)}}\leq\|f\|_{L^{p}(\mathbb{T)}}.
\end{align}
\end{theorem}

Next, we estimate the modulus of each first-order partial derivatives of  $\alpha$-harmonic function $u(z)$ using the $L^{p}$ norm of the boundary function $f$.
\begin{theorem}\label{Th1.7}
Let $u(z)=P_{\alpha}[f](z)$ be $\alpha$-harmonic functions with $f\in L^{p}(\mathbb{T})$, $p\geq 1$ and $z=re^{i\theta}\in\mathbb{D}$. Then we have:

(1)\quad  There exists a function $A_{\alpha, p}(r)$ and  a constant  $A_{\alpha, p}=\sup\limits_{r} A_{\alpha, p}(r)$ such that
\begin{equation}\label{1.13}
|u_{r}|\leq\frac{A_{\alpha, p}(r)}{(1-r^{2})^{1+1/p}}\|f\|_{L^{p}(\mathbb{T)}}\leq\frac{A_{\alpha, p}}{(1-r^{2})^{1+1/p}}\|f\|_{L^{p}(\mathbb{T)}}.
\end{equation}The constant $A_{\alpha, p}$ is asymptotically sharp as $\alpha\rightarrow0$.

(2)\quad  There exists a function $B_{\alpha, p}(r)$ such that \begin{equation}\label{1.14}
|u_{\theta}|\leq\frac{B_{\alpha, p}(r)}{(1-r^{2})^{1+1/p}}\|f\|_{L^{p}(\mathbb{T)}}.
\end{equation}Furthermore, if $\alpha+\frac{2}{p}\geq 0$, then the constant $B_{\alpha, p}=\sup\limits_{r} B_{\alpha, p}(r)$ exists and  is sharp in the inequality
\begin{equation*}
|u_{\theta}|\leq\frac{B_{\alpha, p}}{(1-r^{2})^{1+1/p}}\|f\|_{L^{p}(\mathbb{T)}}.
\end{equation*}

(3)\quad  There exists a function $C_{\alpha, p}(r)$ and a constant  $C_{\alpha, p}=\sup\limits_{r} C_{\alpha, p}(r)$ such that \begin{equation}\label{1.15}
|u_{z}|,  |u_{\bar{z}}|\leq\frac{C_{\alpha, p}(r)}{(1-r^{2})^{1+1/p}}\|f\|_{L^{p}(\mathbb{T)}}\leq\frac{C_{\alpha, p}}{(1-r^{2})^{1+1/p}}\|f\|_{L^{p}(\mathbb{T)}}.
\end{equation}The constant $C_{\alpha, p}$ is asymptotically sharp as $\alpha\rightarrow0$.

 The functions $A_{\alpha, p}(r)$, $B_{\alpha, p}(r)$,  $C_{\alpha, p}(r)$, and the constants  $A_{\alpha, p}$, $B_{\alpha, p}$,  $C_{\alpha, p}$ are defined in (\ref{3.22}), (\ref{3.30}), (\ref{3.37}), (\ref{3.23}), (\ref{3.31}), and (\ref{3.38}), respectively.

\end {theorem}

We can also estimate the $M_{p}(r, \cdot)$ of each first-order partial derivative of  $\alpha$-harmonic function $u(z)$ using the $L^{p}$ norm of the boundary function $f$.

\begin{theorem}\label{Th1.8}
Let $u(z)=P_{\alpha}[f](z)$ be $\alpha$-harmonic functions with $f\in L^{p}(\mathbb{T})$, $p\geq 1$ and $z=re^{i\theta}\in\mathbb{D}$. Then we have:

(1)\quad There exists a function $D_{\alpha}(r)$ and a constant  $D_{\alpha}=\sup\limits_{r} D_{\alpha}(r)$  such that
\begin{align}\label{1.16}
M_{p}(r,u_{r})\leq \frac{D_{\alpha}(r)}{1-r^{2}}\|f\|_{L^{p}(\mathbb{T)}}\leq \frac{D_{\alpha}}{1-r^{2}}\|f\|_{L^{p}(\mathbb{T)}}.
\end{align}The constant $D_{\alpha}$ is asymptotically sharp as $\alpha\rightarrow0$.

(2)\quad There exists a function $E_{\alpha}(r)$ and a constant  $E_{\alpha}=\sup\limits_{r} E_{\alpha}(r)$ such that
\begin{align}\label{1.17}
M_{p}(r,u_{\theta})\leq \frac{E_{\alpha}(r)}{1-r^{2}}\|f\|_{L^{p}(\mathbb{T)}}\leq \frac{E_{\alpha}}{1-r^{2}}\|f\|_{L^{p}(\mathbb{T)}}.
\end{align}The constant $E_{\alpha}$ is asymptotically sharp as $\alpha\rightarrow0$.

(3)\quad There exists a function $F_{\alpha}(r)$ and a constant  $F_{\alpha}=\sup\limits_{r} F_{\alpha}(r)$ such that
\begin{align}\label{1.18}
M_{p}(r,u_{z}), M_{p}(r,u_{\bar{z}})\leq \frac{F_{\alpha}(r)}{1-r^{2}}\|f\|_{L^{p}(\mathbb{T)}}\leq \frac{F_{\alpha}}{1-r^{2}}\|f\|_{L^{p}(\mathbb{T)}}.
\end{align}The constant $F_{\alpha}$ is asymptotically sharp as $\alpha\rightarrow0$.

 The functions $D_{\alpha}(r)$, $E_{\alpha}(r)$, and $F_{\alpha}(r)$, and the constants $D_{\alpha}$, $E_{\alpha}$, and $F_{\alpha}$ are defined in (\ref{3.43}), (\ref{3.49}), (\ref{3.54}), (\ref{3.44}), (\ref{3.50}), and (\ref{3.55}), respectively.

\end{theorem}

{\bf 1.2   To estimate quasiconformal $\alpha$-harmonic functions by the properties of their boundary  functions}

Given $K\geq 1$ and a domain $\Omega$ in $\mathbb{C}$, write $QC(\mathbb{D},\Omega; K)$ for the class of all $K$-quasiconformal mappings  of $\mathbb{D}$ onto $\Omega$ and let $QCH_{\alpha}(\mathbb{D},\Omega; K)$
be the class of all mappings in $QC(\mathbb{D},\Omega; K)$ that  are $\alpha$-harmonic on $\mathbb{D}$.
 when $\alpha=0$, then it is the class of quasiconformal harmonic mappings. As a highly operable type of quasiconformal mapping, harmonic quasi-conformal mapping has received in-depth research and rapid development in recent years, cf.\cite{Kalaj2022,Chuaqui2021,Long2022MM}.

We can use the arc length of the boundary to estimate the integral of the modulus of the first-order partial derivative of the $\alpha$-harmonic quasi-conformal mappings on each concentric circle. This has obvious geometric significance. Specifically, we have

 \begin{theorem}\label{Th1.9}
Let $u\in QCH_{\alpha}(\mathbb{D},\Omega; K)$. If $\Omega$ is bounded by a rectifiable Jordan curve $\gamma$, then
 \begin{equation}\label{1.18}
 \int_{\mathbb{T}_{r}}|\partial u(z)||dz|\leq\frac{K+1}{2}c_{\alpha}F(-\alpha/2, -\alpha/2; 1; r^{2})|\gamma|
\end{equation}and
\begin{equation}\label{1.20}
\int_{\mathbb{T}_{r}}|\bar{\partial} u(z)||dz|\leq\frac{K-1}{2} c_{\alpha}F(-\alpha/2, -\alpha/2; 1; r^{2})|\gamma|,
\end{equation}where $|\gamma|$ is the length of $\gamma$, $\mathbb{T}_{r}=\{z|\,|z|=r\}$.

In particular,
 \begin{equation}\label{1.21}
 \sup_{0<r<1}\int_{\mathbb{T}_{r}}|\partial u(z)||dz|\leq\frac{K+1}{2}|\gamma|
\end{equation}and
\begin{equation}\label{1.22}
\sup_{0<r<1}\int_{\mathbb{T}_{r}}|\bar{\partial} u(z)||dz|\leq\frac{K-1}{2}|\gamma|,
\end{equation}
The constants $\frac{K+1}{2}$ and $\frac{K-1}{2}$ are all sharp when $\alpha\rightarrow 0$ and $K\rightarrow 1$.
\end{theorem}

If the boundary function satisfies the Lipschitz condition, then the modulus of the first-order partial derivative of the $\alpha$-harmonic $K$-quasiconformal mappings can be estimated using Lipschitz constants  $L$ and the dilatation $K$.

 \begin{theorem}\label{Th1.10}
Given $K>1$ and a Jordan domain $\Omega$ in $\mathbb{C}$. Let $u\in QCH_{\alpha}(\mathbb{D}, \Omega, K)$. If the boundary valued function $f$ of $u$ satisfies the inequality
\begin{equation}\label{1.23}
|f(e^{it_{1}})-f(e^{it_{2}})|\leq L |e^{it_{1}}-e^{it_{2}}|,\qquad  t_{1}, t_{2}\in\mathbb{R}
\end{equation}for some positive constant $L$, then
\begin{equation}\label{1.24}
\sup_{z\in\mathbb{D}}|z\partial u(z)|\leq \frac{K+1}{2}L
\end{equation}and
\begin{equation}\label{1.25}
\sup_{z\in\mathbb{D}}|\bar{z}\bar{\partial} u(z)|\leq \frac{K-1}{2}L.
\end{equation}The constants $\frac{K+1}{2}$ and $\frac{K-1}{2}$ are all sharp when $\alpha\rightarrow 0$ and $K\rightarrow 1$.
 \end{theorem}

If $\alpha=0$, then the corresponding results for Theorem \ref{Th1.8} and \ref{Th1.9} can be found in \cite{Partyka2002}.

\section{Prliminaries}
\setcounter{equation}{0}
\hspace{2mm}

In this section, we shall recall some necessary terminology and useful known results. \\

{\bf 2.1 Gauss hypergeometric functions}

The  Gauss hypergeometric function is defined by the  series
$$F(a,b;c; x)=\sum_{n=0}^{\infty}\frac{(a)_{n}(b)_{n}}{(c)_{n}}\frac{x^{n}}{n!}$$ for $ |x|<1$,  and by continuation elsewhere,
where $(a)_{0}=1$ and $(a)_{n}=a(a+1)\cdots(a+n-1)$ for $n=1,2,...$ are the Pochhammer symbols. Obviously, for $n=0, 1,2,...$,
$(a)_{n}=\Gamma(a+n)/\Gamma(a)$.

{\bf  Conclusions}\cite{Andrews1999}

1. If $\Re(c-a-b)>0$, then
\begin{equation}\label{2.1}
\lim_{x\rightarrow 1}F(a,b;c; x)=\frac{\Gamma(c)\Gamma(c-a-b)}{\Gamma(c-a)\Gamma(c-b)}.
\end{equation}

2. If $\Re(c-a-b)<0$ and $|x|<1$, then
\begin{equation}\label{2.2}
F(a,b;c; x)=(1-x)^{c-a-b}F(c-a, c-b; c; x).
\end{equation}

3. It holds that\begin{equation}\label{2.3}
\frac{d F(a,b; c; x)}{dx}=\frac{ab}{c}F(a+1, b+1; c+1; x).
\end{equation}

\begin{lemma}\cite{Olofsson2014}\label{lem2.2}
Let $c>0$, $a\leq c$, $b\leq c$ and $ab\leq 0$ $(ab\geq 0)$. Then the function $F(a,b;c;x)$ is decreasing (increasing) on $x\in (0, 1)$.
\end{lemma}

{\bf 2.2   Jensen's inequality}

Suppose that $\varphi : [\alpha, \beta]\rightarrow\mathbb{R}$ is a convex function, $f$ and $p$ are integrable in $[a,b]$, where $\alpha<\beta$ and $a<b$. If for any $x\in[a,b]$, $f(x)\in[\alpha, \beta]$, $p(x)\geq0$ and $\int^{b}_{a}p(x)dx>0$, then
\begin{equation}\label{2.4}
\varphi\left(\frac{\int^{b}_{a}f(x)p(x)dx}{\int^{b}_{a}p(x)dx}\right)
\leq\frac{\int^{b}_{a}\varphi(f(x))p(x)dx}{\int^{b}_{a}p(x)dx}.
\end{equation}

{\bf 2.3  Hardy space}

Let $f$ be measurable  complex-valued function defined on unit disk $\mathbb{D}$. The integral means of $f$ are defined as follows:
\begin{align*}
M_{p}(r,f)=\left(\frac{1}{2\pi}\int^{2\pi}_{0}|f(re^{i\theta})|^{p}d\theta\right)^{1/p}, \quad 0<p<\infty;
\end{align*}and
\begin{align*}
M_{\infty}(r,f)=ess\sup\limits_{0\leq\theta\leq2\pi}|f(re^{i\theta})|.
\end{align*}
A function $f$ analytic in $\mathbb{D}$ is said to be of class $H^{p}(\mathbb{D})$ if $M_{p}(r,u)$ is bounded.

It is convenient also to define the analogous classes of harmonic functions or $\alpha$-harmonic functions. A function $u$ harmonic or $\alpha$-harmonic in $\mathbb{D}$ is said to be of class $h^{p}(\mathbb{D})$ or $h_{\alpha}^{p}(\mathbb{D})$ if $M_{p}(r,u)$ is bounded, respectively.\\

{\bf 2.4  Quasiconformal mappings}

\begin{definition} Suppose $f:  \Omega\rightarrow\Omega'$ is a homeomorphic $W_{loc}^{1,2}$-mapping. Then $f$ is $K$-quasiconformal if
\begin{align*}
f_{\bar{z}}(z)=\mu(z)f_{z}(z)
\end{align*}for almost every $z\in\Omega$, where $\mu$, call the Beltrami coefficient of $f$, is a bounded measurable function satisfying
\begin{align*}\|\mu\|\leq\frac{K-1}{K+1}<1.
\end{align*}
\end{definition}

{\bf 2.5  An important  variable substitution}

 The following variable substitution should be used in the proof of Theorem \ref{Th1.6} and \ref{Th1.7} for several times.

For $t, \theta, s\in[0,2\pi]$, let \begin{align*}
e^{i(t-\theta)}=\frac{r+e^{is}}{1+re^{is}}.
\end{align*}Then we have
 \begin{equation*}
|1-re^{i(\theta-t)}|=\frac{1-r^{2}}{|1+re^{-is}|},
\end{equation*}
\begin{equation*}
r-\cos(\theta-t)=r-\frac{2r+(1+r^{2})\cos s}{|1+re^{-is}|^{2}}=\frac{-(1-r^{2})(r+\cos s)}{|1+re^{-is}|^{2}},
\end{equation*}and
\begin{align*}
dt=\frac{1-r^{2}}{|1+re^{is}|^{2}}ds.
\end{align*}

{\bf 2.6  Several Lemmas}

\begin{lemma}\cite{Kalaj2024} \label{lem2.3}
Let $m>-1$ and $k\geq 0$ and define
\begin{align}G(r,\theta):=\int^{\pi}_{-\pi}|\cos(t-\theta)|^{k}(1+r^{2}+2r\cos t)^{m}dt.
\end{align}Then \begin{equation} \; G(r,\theta)\leq \left\{
\begin{array}{rr}
G(1,0),  &\quad m>1,\\
G(1,\frac{\pi}{2}),  &\quad m\leq 1.
\end{array}\right.
\end{equation}
\end{lemma}

\begin{lemma}(P.1179 of \cite{Kalaj2013})\label{lem2.4}

It holds that
\begin{align}\label{2.7}
\int^{\pi}_{0}\frac{\sin^{\mu-1}t}{(1+r^{2}-2r\cos t)^{\nu}}dt=B(\frac{\mu}{2},\frac{1}{2})F(\nu,  \nu+\frac{1-\mu}{2}; \frac{1+\mu}{2}; r^{2}),
\end{align}where $B(u,v)$ is the Beta-function, and $F(a, b; c; x)$ is the hypergeometric Gauss function.
Specially,
\begin{align}\label{2.8}
\int^{\pi}_{0}\sin^{\mu-1}t\,(1-\cos t)^{-\nu}dt=2^{\nu}B(\frac{\mu}{2},\frac{1}{2})F(\nu,  \nu+\frac{1-\mu}{2}; \frac{1+\mu}{2}; 1),
\end{align}
and
\begin{align}\label{2.9}
\int^{\pi}_{0}\frac{dt}{(1+r^{2}-2r\cos t)^{\nu}}=\pi F(\nu,  \nu; 1; r^{2}).
\end{align}
\end{lemma}

\section{Proofs}
\setcounter{equation}{0}
\hspace{2mm}

\begin{proof}[\textbf{Proof of Theorem \ref{Th1.6}}]

Rewrite equation (\ref{1.3}) to
\begin{align}\label{3.15}
u(z)=\frac{1}{2\pi}\int^{2\pi}_{0}c_{\alpha}\frac{(1-r^{2})^{\alpha+1}}{|1-re^{i(\theta-t)}|^{\alpha+2}}f(e^{it})dt
,\quad z=re^{i\theta}. \end{align}
Then we have \begin{align}\label{3.16}
|u(z)|\leq \int^{2\pi}_{0}c_{\alpha}\frac{(1-r^{2})^{\alpha+1}}{|1-re^{i(\theta-t)}|^{\alpha+2}}|f(e^{it})|\frac{dt}{2\pi}
. \end{align}
Let \begin{align*}
I= \int^{2\pi}_{0}c_{\alpha}\frac{(1-r^{2})^{\alpha+1}}{|1-re^{i(\theta-t)}|^{\alpha+2}}\frac{dt}{2\pi}
. \end{align*}
Then from the proof process of the Theorem 3.1 of \cite{Olofsson2014}, we obtain that
 \begin{align*}
I=c_{\alpha}F(-\alpha/2, -\alpha/2; 1; r^{2}).
\end{align*}By Lemma \ref{lem2.2} and equation (\ref{2.1}), we have
 \begin{align*}
I=c_{\alpha}F(-\alpha/2, -\alpha/2; 1; r^{2})\leq 1.
\end{align*}
For $p\geq 1$, considering Jensen's inequality for (\ref{3.16}), we have
\begin{align*}
|u(z)|^{p}&\leq \left(I\cdot \int^{2\pi}_{0}\frac{c_{\alpha}\frac{(1-r^{2})^{\alpha+1}}{|1-re^{i(\theta-t)}|^{\alpha+2}}}{I}|f(e^{it})|\frac{dt}{2\pi}\right)^{p}\\
&\leq I^{p-1}\int^{2\pi}_{0}c_{\alpha}\frac{(1-r^{2})^{\alpha+1}}{|1-re^{i(\theta-t)}|^{\alpha+2}}|f(e^{it})|^{p}\frac{dt}{2\pi}
. \end{align*}
It follows that
\begin{align*}
\frac{1}{2\pi}\int^{2\pi}_{0}|u(z)|^{p}d\theta &\leq I^{p-1}\frac{1}{2\pi}\int^{2\pi}_{0}\left(\int^{2\pi}_{0}c_{\alpha}\frac{(1-r^{2})^{\alpha+1}}{|1-re^{i(\theta-t)}|^{\alpha+2}}|f(e^{it})|^{p}\frac{dt}{2\pi}\right)d\theta\\
&\leq I^{p}\|f\|^{p}_{L^{p}(\mathbb{T})}.
\end{align*}
Therefore,\begin{align*} M_{p}(r,u)\leq I\cdot\|f\|_{L^{p}(\mathbb{T})}= c_{\alpha}F(-\alpha/2, -\alpha/2; 1; r^{2})\|f\|_{L^{p}(\mathbb{T)}}\leq\|f\|_{L^{p}(\mathbb{T)}}.\end{align*}

Let $f(e^{it})\equiv C$ with $C>0$. Then (\ref{3.15}) shows that $u(z)= C\cdot c_{\alpha}F(-\alpha/2, -\alpha/2; 1; r^{2})$. Thus it is easy to see that both inequalities (\ref{1.11}) and (\ref{1.12}) are sharp.
\end{proof}

\begin{proof} [\textbf{Proof of Theorem \ref{Th1.7}}]

(1)   Differential on both sides of the formula
(\ref{3.15}) , we have
\begin{align}\label{3.17}
\frac{\partial}{\partial r}u(re^{i\theta})
=-\frac{c_{\alpha}}{2\pi}(1-r^{2})^{\alpha}\int^{2\pi}_{0}\frac{2(\alpha+1)r|1-re^{i(\theta-t)}|^{2}+(\alpha+2)(1-r^{2})(r-\cos(\theta-t))}
{|1-re^{i(\theta-t)}|^{4+\alpha}}f(e^{it})dt.
\end{align}
Let \begin{align*}
I_{1}=\int^{2\pi}_{0}\frac{|2(\alpha+1)r|1-re^{i(\theta-t)}|^{2}+(\alpha+2)(1-r^{2})(r-\cos(\theta-t))|^{q}}
{|1-re^{i(\theta-t)}|^{(4+\alpha)q}}\frac{dt}{2\pi}
\end{align*}with $\frac{1}{p}+\frac{1}{q}=1$.
Then by H\"{o}lder inequality and equation (\ref{3.17}) we have
\begin{align}\label{3.18}
|u_{r}|\leq c_{\alpha}(1-r^{2})^{\alpha}I_{1}^{\frac{1}{q}}\|f\|_{L^{p}(\mathbb{T)}}.
\end{align}
Make variable substitution as in section 2.5, we have
\begin{align}\label{3.19}
I_{1}=&\int^{2\pi}_{0}\frac{\left(\frac{1-r^{2}}{|1+re^{-is}|}\right)^{2q}|2(\alpha+1)r-(\alpha+2)(r+\cos s)|^{q}}{\left(\frac{1-r^{2}}{|1+re^{-is}|}\right)^{(4+\alpha)q}}
\frac{1-r^{2}}{|1+re^{is}|^{2}}\frac{ds}{2\pi}\nonumber\\
=&(1-r^{2})^{1-(2+\alpha)q}\int^{2\pi}_{0}\frac{|\alpha r-(\alpha+2)\cos s|^{q}}{|1+re^{-is}|^{2-(2+\alpha)q}}\frac{ds}{2\pi}.
\end{align}
Let
\begin{align}\label{3.20}
\widetilde{A}_{\alpha, p}(r)=c_{\alpha}\left(\int^{2\pi}_{0}\frac{|\alpha r-(\alpha+2)\cos s|^{q}}{|1+re^{-is}|^{2-(2+\alpha)q}}\frac{ds}{2\pi}\right)^{1/q}.
\end{align}
Then
(\ref{3.18}) and (\ref{3.19}) leads to
\begin{align}\label{3.21}
|u_{r}|\leq \frac{\widetilde{A}_{\alpha, p}(r)}{(1-r^{2})^{1+1/p}}\|f\|_{L^{p}(\mathbb{T)}}.
\end{align}

Recall the formula of Lagrange mean value theorem $f(x_{0}+ x)=f(x_{0})+f'(\xi)x$. It follows that $f(x_{0}+ x)\leq f(x_{0})+ |x|\max\limits_{\xi}f'(\xi)$. Then we have
\begin{align*}
&\quad|\alpha r-(\alpha+2)\cos s|^{q}\leq (|\alpha| r+|(\alpha+2)\cos s|)^{q}\\
&\leq |(\alpha+2)\cos s|^{q}+\max\{q(|\alpha| r+|(\alpha+2)\cos s|)^{q-1}\}|\alpha| r\\
&\leq (\alpha+2)^{q}|\cos s|^{q}+q(|\alpha| r+\alpha+2)^{q-1}|\alpha| r.
\end{align*}
So for the term of (\ref{3.20}),  we have the following estimate:
\begin{align*}
&\quad \int^{2\pi}_{0}\frac{|\alpha r-(\alpha+2)\cos s|^{q}}{|1+re^{-is}|^{2-(2+\alpha)q}}\frac{ds}{2\pi}\\
&\leq \frac{ (\alpha+2)^{q}}{2\pi}\int^{2\pi}_{0}\frac{|\cos s|^{q}}{|1+re^{-is}|^{2-(2+\alpha)q}}ds+\frac{q(|\alpha| r+\alpha+2)^{q-1}|\alpha| r}{2\pi}\int^{2\pi}_{0}\frac{1}{|1+re^{-is}|^{2-(2+\alpha)q}}ds\\
&:=\frac{ (\alpha+2)^{q}}{2\pi}I_{11}+\frac{q(|\alpha| r+\alpha+2)^{q-1}|\alpha| r}{2\pi}I_{12}.
\end{align*}

Rewrite $I_{11}$ as
\begin{align*}
&I_{11}=\int^{2\pi}_{0}|\cos s|^{q}|1+re^{-is}|^{(2+\alpha)q-2}ds=\int^{2\pi}_{0}|\cos s|^{q}(1+r^{2}+2r\cos s)^{\frac{(2+\alpha)q-2}{2}}ds=:G(r,0).
\end{align*}
Now we  discuss it as follows:

If $(2+\alpha)q>4$, then $\frac{(2+\alpha)q-2}{2}>1$. By Lemma \ref{lem2.3}, we have
$$G(r,0)\leq G(1,0)=\int^{2\pi}_{0}|\cos s|^{q}(2+2\cos s)^{(2+\alpha)q-2}ds.$$

If  $(2+\alpha)q\leq4$, then $\frac{(2+\alpha)q-2}{2}\leq1$. By Lemma \ref{lem2.3}, we have
\begin{align*}
G(r,0)&\leq G(1,\frac{\pi}{2})=\int^{2\pi}_{0}|\cos (s-\frac{\pi}{2})|^{q}(2+2\cos s)^{(2+\alpha)q-2}ds\\
&=\int^{2\pi}_{0}|\sin s|^{q}(2+2\cos s)^{(2+\alpha)q-2}ds\\
&=\int^{\pi}_{-\pi}|\sin t|^{q}(2-2\cos t)^{(2+\alpha)q-2}dt\\
&=2^{(2+\alpha)q-1}\int^{\pi}_{0}(\sin t)^{q}(1-\cos t)^{(2+\alpha)q-2}dt.
\end{align*}
Consider (\ref{2.8}) of Lemma \ref{lem2.4}, we have
\begin{align*}G(1,\frac{\pi}{2})
&=2 B\left(\frac{q+1}{2},\frac{1}{2}\right)F\left(2-(2+\alpha)q, 2-\frac{(5+2\alpha)q}{2}; 1+\frac{q}{2}; 1\right)=\frac{2\sqrt{\pi} \Gamma(\frac{q+1}{2})\Gamma((2\alpha+5)q-3)}{\Gamma(\frac{(\alpha+5)q}{2}-1)\Gamma((\alpha+3)q-1)}.
\end{align*}Here, the second equality holds just because (\ref{2.1}) holds  and $2-(2+\alpha)q+2-\frac{(5+2\alpha)q}{2}< 1+\frac{q}{2}$ holds  for $q\geq1$ and $\alpha>-1$.

Rewrite $I_{12}$ as
\begin{align*}
&I_{12}=\int^{2\pi}_{0}|1+re^{-is}|^{(2+\alpha)q-2}ds=\int^{\pi}_{-\pi}(1+r^{2}-2r\cos s)^{\frac{(2+\alpha)q-2}{2}}dt=2\int^{\pi}_{0}(1+r^{2}-2r\cos s)^{\frac{(2+\alpha)q-2}{2}}dt.
\end{align*}
Considering (\ref{2.9}) of Lemma \ref{lem2.4}, we have
\begin{align*}
I_{12}&=2\pi F\left(\frac{2-(2+\alpha)q}{2}, \frac{2-(2+\alpha)q}{2}; 1; r^{2}\right)
< 2\pi F\left(\frac{2-(2+\alpha)q}{2}, \frac{2-(2+\alpha)q}{2}; 1; 1\right)
=2\pi\frac{\Gamma((\alpha+2)q-1)}{\left(\Gamma(\frac{(\alpha+2)q}{2})\right)^{2}}.
\end{align*}Here, the inequality holds  because of Lemma \ref{lem2.2}; the second equality holds just because (\ref{2.1}) holds  and $\frac{2-(2+\alpha)q}{2}+\frac{2-(2+\alpha)q}{2}< 1$ holds  for $q\geq1$ and $\alpha>-1$.

Let \begin{equation}\label{3.22}
A_{\alpha, p}(r)=c_{\alpha}\left(\frac{ (\alpha+2)^{q}}{2\pi}I_{11}+\frac{q(|\alpha| r+\alpha+2)^{q-1}|\alpha| r}{2\pi}I_{12}\right)^{\frac{1}{q}}
\end{equation}and
\begin{equation} \label{3.23}\; A_{\alpha, p}= \left\{
\begin{array}{rr}
c_{\alpha}\left(\frac{ (\alpha+2)^{q}}{2\pi}G(1,0)+q(|\alpha|+\alpha+2)^{q-1}|\alpha| \frac{\Gamma((\alpha+2)q-1)}{\left(\Gamma(\frac{(\alpha+2)q}{2})\right)^{2}}\right)^{\frac{1}{q}}, &\quad (2+\alpha)q>4;\\
c_{\alpha}\left(\frac{ (\alpha+2)^{q}}{2\pi}G(1,\frac{\pi}{2})+q(|\alpha|+\alpha+2)^{q-1}|\alpha| \frac{\Gamma((\alpha+2)q-1)}{\left(\Gamma(\frac{(\alpha+2)q}{2})\right)^{2}}\right)^{\frac{1}{q}},  &\quad (2+\alpha)q\leq 4.
\end{array}\right.
\end{equation}
Then we have    $\tilde{A}_{\alpha, p}(r)\leq A_{\alpha, p}(r)\leq A_{\alpha, p}$. Therefore, it follows from   (\ref{3.21}) that  (\ref{1.13}) holds with $A_{\alpha, p}(r)\leq A_{\alpha, p}$.

Specially, if $\alpha=0$, then
\begin{align*} A_{0, p}(r)=2\left(\frac{ 1}{2\pi}\int^{2\pi}_{0}|\cos s|^{q}|1+re^{-is}|^{2q-2}ds\right)^{\frac{1}{q}},\\
 A_{0, p}=\frac{4^{\frac{1}{p}}}{\pi^{\frac{1}{q}}}\left(\int^{2\pi}_{0}|\cos s|^{q}(1+\cos s)^{q-1}ds\right)^{\frac{1}{q}}.
  \end{align*}

Next, let us show the constant $A_{0, p}$  is sharp.  Let $0<\rho<1$, and define
\begin{align}\label{3.24}
f_{\rho}(e^{it})=(1-\rho^{2})^{-\frac{1}{p}}|\cos s|^{q-1}(1+\cos s)^{q-1}sign(\cos s),
\end{align}where $s$ and $t$ satisfy the relation
\begin{align}\label{3.25}
e^{i(t-\theta)}=\frac{\rho+e^{is}}{1+\rho e^{is}},
\end{align}which is similar to the variable substitution of section 2.5.
Then we have
\begin{align*}\frac{1}{2\pi}\int^{2\pi}_{0}|f_{\rho}(e^{it})|^{p}dt
&=\frac{1}{2\pi}\int^{2\pi}_{0}(1-\rho^{2})^{-1}|\cos s|^{p(q-1)}(1+\cos s)^{p(q-1)}\frac{1-\rho^{2}}{|1+\rho e^{is}|^{2}}ds\\
&=\frac{1}{2\pi}\int^{2\pi}_{0}|\cos s|^{q}(1+\cos s)^{q}|1+\rho e^{is}|^{-2}ds.
\end{align*}
Thus, \begin{align*}&\lim_{\rho\rightarrow 1}\|f_{\rho}\|^{p}_{p}=\lim_{\rho\rightarrow 1}
\frac{1}{2\pi}\int^{2\pi}_{0}|\cos s|^{q}(1+\cos s)^{q}|1+\rho e^{is}|^{-2}ds
=\frac{1}{4\pi}\int^{2\pi}_{0}|\cos s|^{q}(1+\cos s)^{q-1}ds.
\end{align*}
Now, take $$u(z)=P_{\alpha}[f_{\rho}](z).$$
If $\alpha=0$, then (\ref{3.17}) reduces to
\begin{align*}
\frac{\partial}{\partial r}u(re^{i\theta})
&=-\frac{1}{\pi}\int^{2\pi}_{0}\frac{r|1-re^{i(\theta-t)}|^{2}+(1-r^{2})(r-\cos(\theta-t))}
{|1-re^{i(\theta-t)}|^{4}}f(e^{it})dt.
\end{align*}
Taking $r=\rho$, and making variable change as (\ref{3.25}), we have
\begin{align*}
\frac{\partial}{\partial r}u(re^{i\theta})
&=-\frac{1}{\pi}\int^{2\pi}_{0}\frac{\rho|1-\rho\frac{\rho+e^{-is}}{1+\rho e^{-is}}|^{2}-(1-\rho^{2})\frac{(1-\rho^{2})(\rho+\cos s)}{|1+\rho e^{-is}|^{2}}}
{|1-\rho\frac{\rho+e^{-is}}{1+\rho e^{-is}}|^{4}}f(e^{it})\frac{1-\rho^{2}}{|1+\rho e^{-is}|^{2}}ds\\
&=\frac{1}{\pi}\int^{2\pi}_{0}\frac{\cos s}{1-\rho^{2}}f(e^{it})ds\\
&=\frac{1}{\pi}\int^{2\pi}_{0}\frac{\cos s}{1-\rho^{2}}(1-\rho^{2})^{-\frac{1}{p}}|\cos s|^{q-1}(1+\cos s)^{q-1}sign(\cos s)ds\\
&=(1-\rho^{2})^{-\frac{1}{p}-1}\frac{1}{\pi}\int^{2\pi}_{0}|\cos s|^{q}(1+\cos s)^{q-1}ds.
\end{align*}
Therefore,
\begin{align*}
\lim_{r=\rho\rightarrow 1}\frac{(1-\rho^{2})^{1+\frac{1}{p}}|\frac{\partial}{\partial r}u|}{\|f_{\rho}\|_{p}}=A_{0,p}.
\end{align*}
This shows that the constant $A_{0,p}$ is sharp.

(2)    Differential on both sides of the formula
(\ref{3.15}) , we have
\begin{align}\label{3.26}
\frac{\partial}{\partial \theta}u(re^{i\theta})&=-\frac{(\alpha+2)c_{\alpha}}{2\pi}r(1-r^{2})^{\alpha+1}\int^{2\pi}_{0}\frac{\sin(\theta-t)}
{(1+r^{2}-2r\cos(\theta-t))^{2+\alpha/2}}f(e^{it})dt.
\end{align}
Let
\begin{align}\label{3.27}
I_{2}=\int^{2\pi}_{0}\frac{|\sin(\theta-t)|^{q}}
{(1+r^{2}-2r\cos(\theta-t))^{(2+\alpha/2)q}}\frac{dt}{2\pi}
\end{align}with $\frac{1}{p}+\frac{1}{q}=1$.
Then by H\"{o}lder inequality, from (\ref{3.26}), we have
\begin{align}\label{3.28}
|u_{\theta}|\leq (\alpha+2)c_{\alpha}r(1-r^{2})^{\alpha+1}I_{2}^{\frac{1}{q}}\|f\|_{L^{p}(\mathbb{T)}}.
\end{align}
Using
 (\ref{2.7}),  equation (\ref{3.27}) reduces to
  \begin{align}\label{3.29}
I_{2}&=\int^{2\pi}_{0}\frac{|\sin t|^{q}}
{(1+r^{2}-2r\cos t)^{(2+\alpha/2)q}}\frac{dt}{2\pi}\nonumber\\
&=\frac{1}{\pi}\int^{\pi}_{0}\frac{\sin^{q} t}
{(1+r^{2}-2r\cos t)^{(2+\alpha/2)q}}dt\nonumber\\
&=\frac{1}{\pi}B\left(\frac{q+1}{2},\frac{1}{2}\right)F\left((\frac{\alpha}{2}+2)q,\frac{(\alpha+3)q}{2}; 1+\frac{q}{2};r^{2}\right)\nonumber\\
&=\frac{1}{\pi}B\left(\frac{q+1}{2},\frac{1}{2}\right)F\left((1-\frac{(\alpha+3)q}{2},1-\frac{(\alpha+2)q}{2}; 1+\frac{q}{2};r^{2}\right)(1-r^{2})^{1-(\alpha+3)q}.
\end{align}Here, the last equality comes from (\ref{2.2}) because $(\frac{\alpha}{2}+2)q+\frac{(\alpha+3)q}{2}>1+\frac{q}{2}$ for all $\alpha>-1$ and $p\geq 1$.
Let \begin{align}\label{3.30}
B_{\alpha, p}(r)=\frac{(\alpha+2)c_{\alpha}}{\pi^{1/q}}r\left(B\left(\frac{q+1}{2},\frac{1}{2}\right)
F\left((1-\frac{(\alpha+3)q}{2},1-\frac{(\alpha+2)q}{2}, 1+\frac{q}{2};r^{2}\right)\right)^{1/q}.
\end{align}Then
(\ref{3.28}) and (\ref{3.29}) lead to
\begin{align*}
|u_{\theta}|\leq \frac{B_{\alpha, p}(r)}{(1-r^{2})^{1+1/p}}\|f\|_{L^{p}(\mathbb{T)}}.
\end{align*}

If $\alpha+\frac{2}{p}\geq0$, then $1-\frac{(\alpha+3)q}{2}+1-\frac{(\alpha+2)q}{2}<1+\frac{q}{2}$  and $\left(1-\frac{(\alpha+3)q}{2}\right)\left(1-\frac{(\alpha+2)q}{2}\right)\geq0$. By Lemma \ref{lem2.2}, it holds that
\begin{align}\label{3.31}
B_{\alpha, p}(r)&\leq\frac{(\alpha+2)c_{\alpha}}{\pi^{1/q}}\left(B\left(\frac{q+1}{2},\frac{1}{2}\right)
F\left(1-\frac{(\alpha+3)q}{2},1-\frac{(\alpha+2)q}{2}; 1+\frac{q}{2};1\right)\right)^{1/q}:=B_{\alpha, p}.
\end{align}Thus, it holds that
\begin{align*}
|u_{\theta}|\leq \frac{B_{\alpha, p}}{(1-r^{2})^{1+1/p}}\|f\|_{L^{p}(\mathbb{T)}}.
\end{align*}

Now we show the constant $B_{\alpha, p}$ is sharp in the above inequality.
Let $0<\rho<1$, and define
\begin{align}
f_{\rho}(e^{it})=(1-\rho^{2})^{\frac{3}{1-p}}|\sin s|^{\frac{1}{p-1}}(1+\cos s)^{\frac{\frac{\alpha}{2}+1}{p-1}}sign(\sin s),
\end{align}where $s$ and $t$ satisfy the relation (\ref{3.25}).
Then \begin{align*}\frac{1}{2\pi}\int^{2\pi}_{0}|f_{\rho}(e^{it})|^{p}dt
&=\frac{1}{2\pi}\int^{2\pi}_{0}(1-\rho^{2})^{\frac{3p}{1-p}}|\sin s|^{\frac{p}{p-1}}(1+\cos s)^{\frac{(\frac{\alpha}{2}+1)p}{p-1}}\frac{1-\rho^{2}}{|1+\rho e^{is}|^{2}}ds\\
&=\frac{1}{2\pi}(1-\rho^{2})^{\frac{2p+1}{1-p}}\int^{2\pi}_{0}|\sin s|^{\frac{p}{p-1}}(1+\cos s)^{\frac{(\frac{\alpha}{2}+1)p}{p-1}}|1+\rho e^{is}|^{-2}ds.
\end{align*}

 Take $$u(z)=P_{\alpha}[f_{\rho}](z)$$with  $\alpha+\frac{2}{p}\geq0$.
Equation (\ref{3.26}) and variable  substitution (\ref{3.25}) lead to
\begin{align*}
\frac{\partial}{\partial \theta}u(re^{i\theta})&=-\frac{(\alpha+2)c_{\alpha}}{2\pi}r(1-r^{2})^{\alpha+1}\int^{2\pi}_{0}\frac{\sin(\theta-t)}
{|1-re^{i(\theta-t)}|^{4+\alpha}}f(e^{it})dt\\
&=-\frac{(\alpha+2)c_{\alpha}}{2\pi}r(1-r^{2})^{\alpha+1}\int^{2\pi}_{0}
\frac{1}{2i}\left(\frac{r+e^{-is}}{1+re^{-is}}-\frac{r+e^{is}}{1+re^{is}}\right)
\frac{|1+re^{-is}|^{4+\alpha}}{(1-r^{2})^{4+\alpha}}f(e^{it})\frac{1-r^{2}}{|1+re^{is}|^{2}}ds\\
&=\frac{(\alpha+2)c_{\alpha}}{2\pi}r(1-r^{2})^{-1}\int^{2\pi}_{0}\sin s \,|1+re^{is}|^{\alpha}f(e^{it})ds\\
&=\frac{(\alpha+2)c_{\alpha}}{2\pi}r(1-r^{2})^{-1}(1-\rho^{2})^{\frac{3}{1-p}}\int^{2\pi}_{0}|\sin s|^{\frac{p}{p-1}}(1+\cos s)^{\frac{\frac{\alpha}{2}+1}{p-1}}|1+re^{is}|^{\alpha}ds.
\end{align*}

Let $r=\rho\rightarrow 1$, then we have
 \begin{align*}
&\quad\lim_{r=\rho\rightarrow 1}\frac{(1-r^{2})^{1+\frac{1}{p}}\,u_{\theta}}{\|f_{\rho}\|_{p}}\\
&=2^{\frac{1}{p}-\frac{1}{q}+\frac{\alpha}{2}}\pi^{-\frac{1}{q}}(\alpha+2)c_{\alpha}\left(\int^{2\pi}_{0}|\sin s|^{\frac{p}{p-1}}(1+\cos s)^{\frac{\frac{\alpha}{2}p+1}{p-1}}ds\right)^{\frac{1}{q}}\\
&=2^{\frac{1}{p}+\frac{\alpha}{2}}\pi^{-\frac{1}{q}}(\alpha+2)c_{\alpha}\left(\int^{\pi}_{0}(\sin s)^{\frac{p}{p-1}}(1-\cos s)^{\frac{\frac{\alpha}{2}p+1}{p-1}}ds\right)^{\frac{1}{q}}\\
&=2^{\frac{1}{p}+\frac{\alpha}{2}}\pi^{-\frac{1}{q}}(\alpha+2)c_{\alpha}
\left(2^{\frac{\frac{\alpha}{2}p+1}{1-p}}B\left(\frac{2p-1}{2(p-1)},\frac{1}{2}\right)
F\left(\frac{\frac{\alpha}{2}p+1}{1-p},\frac{(\alpha+1)p+2}{2(1-p)}; 1+\frac{p}{2(p-1)};1\right)\right)^{\frac{1}{q}}\\
&=\frac{(\alpha+2)c_{\alpha}}{\pi^{1/q}}\left(B\left(\frac{q+1}{2},\frac{1}{2}\right)
F\left(1-\frac{(\alpha+3)q}{2},1-\frac{(\alpha+2)q}{2}, 1+\frac{q}{2};1\right)\right)^{1/q}\\
&=B_{\alpha,p}.
\end{align*}Here the third equality holds  because of (\ref{2.8}). This means the constant $B_{\alpha,p}$  is sharp  when $\alpha+\frac{2}{p}\geq0$.

(3) Differential on both sides of the formula
(\ref{3.15}), we have
\begin{align}\label{3.33}
\frac{\partial}{\partial z}u(z)=c_{\alpha}(1-r^{2})^{\alpha}\int^{2\pi}_{0}\frac{-(\alpha+1)\bar{z}|1-ze^{-it}|^{2}
+(1+\alpha/2)(1-r^{2})e^{-it}(1-\bar{z}e^{it})}{|1-ze^{-it}|^{\alpha+4}}
f(e^{it})\frac{dt}{2\pi}.
\end{align}
Let \begin{align}\label{3.34}
I_{3}&=\int^{2\pi}_{0}\frac{|-(\alpha+1)\bar{z}|1-ze^{-it}|^{2}
+(1+\alpha/2)(1-r^{2})e^{-it}(1-\bar{z}e^{it})|^{q}}{|1-ze^{-it}|^{(\alpha+4)q}}
\frac{dt}{2\pi}\nonumber\\
&=\int^{2\pi}_{0}\frac{|-(\alpha+1)\bar{z}e^{it}|1-ze^{-it}|^{2}
+(1+\alpha/2)(1-r^{2})(1-\bar{z}e^{it})|^{q}}{|1-ze^{-it}|^{(\alpha+4)q}}
\frac{dt}{2\pi}
\end{align}with $\frac{1}{p}+\frac{1}{q}=1$.
Then by H\"{o}lder inequality, we have
\begin{align}\label{3.35}
|u_{z}|\leq  c_{\alpha}(1-r^{2})^{\alpha}I_{3}^{\frac{1}{q}}\|f\|_{L^{p}(\mathbb{T)}}.
\end{align}
After the change of variables same as in section 2.5,  (\ref{3.34}) reduces to
 \begin{align*}
 I_{3}&=\int^{2\pi}_{0}\frac{|-(\alpha+1)\bar{z}e^{it}(\frac{1-r^{2}}{|1+re^{-is}|})^{2}
+(1+\alpha/2)(1-r^{2})(1-r\frac{r+e^{is}}{1+re^{is}})|^{q}}{(\frac{1-r^{2}}{|1+re^{-is}|})^{(\alpha+4)q}}\frac{1-r^{2}}{|1+re^{is}|^{2}}
\frac{ds}{2\pi}\nonumber\\
&=(1-r^{2})^{1-(\alpha+2)q}\int^{2\pi}_{0}\frac{|-(\alpha+1)re^{i(t-\theta)}
+(1+\alpha/2)(r+e^{-is})e^{i(t-\theta)}|^{q}}{|1+re^{-is}|^{2-(\alpha+2)q}}\frac{ds}{2\pi}\nonumber\\
&=(1-r^{2})^{1-(\alpha+2)q}\int^{2\pi}_{0}\frac{|-\alpha r/2
+(1+\alpha/2)e^{-is}|^{q}}{|1+re^{-is}|^{2-(\alpha+2)q}}\frac{ds}{2\pi}.
\end{align*}
Let
 \begin{equation*}
\widetilde{C}_{\alpha, p}(r)=c_{\alpha}\left(\int^{2\pi}_{0}\frac{|-\alpha r/2
+(1+\alpha/2)e^{-is}|^{q}}{|1+re^{-is}|^{2-(\alpha+2)q}}
\frac{ds}{2\pi}\right)^\frac{1}{q}.
\end{equation*}
Then (\ref{3.35})  leads to
\begin{align}\label{3.36}
|u_{z}|\leq \frac{\widetilde{C}_{\alpha, p}(r)}{(1-r^{2})^{1+1/p}}\|f\|_{L^{p}(\mathbb{T)}}.
\end{align}
Observe that $|-\alpha r/2+(1+\alpha/2)e^{-is}|^{q}\leq (|\alpha| r/2+1+\alpha/2)^{q}$. It follows that
\begin{equation}\label{3.37}
\widetilde{C}_{\alpha, p}(r)\leq c_{\alpha}(|\alpha |r/2
+1+\alpha/2)\left(\frac{1}{2\pi}\int^{2\pi}_{0}\frac{ds}{|1+re^{-is}|^{2-(\alpha+2)q}}
\right)^{\frac{1}{q}}=:C_{\alpha, p}(r).
\end{equation}
Considering (\ref{2.9}), we rewrite $C_{\alpha, p}(r)$ as
\begin{align*}
C_{\alpha, p}(r)&=c_{\alpha}(|\alpha|r/2+1+\alpha/2)
\left(\frac{1}{\pi}\int^{\pi}_{0}\frac{dT}{|1-re^{-iT}|^{2-(\alpha+2)q}}\right)^{\frac{1}{q}}\\
&=c_{\alpha}(|\alpha |r/2+1+\alpha/2)
\left( F(1-\frac{(\alpha+2)q}{2}, 1-\frac{(\alpha+2)q}{2}; 1; r^{2})\right)^{\frac{1}{q}}.
\end{align*}
Let \begin{align}\label{3.38}
C_{\alpha, p}&:=c_{\alpha}(|\alpha |r/2+1+\alpha/2)
\left( F(1-\frac{(\alpha+2)q}{2}, 1-\frac{(\alpha+2)q}{2}; 1; 1)\right)^{\frac{1}{q}}\nonumber\\
&= c_{\alpha}(|\alpha |/2
+1+\alpha/2)\left(\frac{\Gamma((\alpha+2)q-1)}{\left(\Gamma(\frac{(\alpha+2)q}{2})\right)^{2}}
\right)^\frac{1}{q}.
\end{align} Then  (\ref{1.15}) follows from (\ref{3.36}) and  (\ref{3.37}) with $C_{\alpha, p}(r)\leq C_{\alpha, p}$.

If $\alpha=0$, then $C_{0, p}(r)=\left(\frac{1}{2\pi}\int^{2\pi}_{0}\frac{ds}{|1+re^{-is}|^{2-2q}}
\right)^\frac{1}{q}$ and $C_{0, p}=\left(\frac{\Gamma(2q-1)}{\left(\Gamma(q)\right)^{2}}
\right)^\frac{1}{q}$. Next, we show the constant $C_{0,p}$ is sharp.

 Let $0<\rho<1$, and define
\begin{align}
f_{\rho}(e^{it})=(1-\rho^{2})^{\frac{2}{1-p}}(1+\cos s)^{\frac{1}{p-1}}e^{is},
\end{align}where $s$ and $t$ satisfy the relation (\ref{3.25}).
Then direct computation leads to
\begin{align*}\frac{1}{2\pi}\int^{2\pi}_{0}|f_{\rho}(e^{it})|^{p}dt
&=\frac{1}{2\pi}\int^{2\pi}_{0}(1-\rho^{2})^{\frac{2p}{1-p}}(1+\cos s)^{\frac{p}{p-1}}\frac{1-\rho^{2}}{|1+\rho e^{is}|^{2}}ds\\
&=\frac{1}{2\pi}(1-\rho^{2})^{\frac{p+1}{1-p}}\int^{2\pi}_{0}(1+\cos s)^{\frac{p}{p-1}}|1+\rho e^{is}|^{-2}ds.
\end{align*}

 Take $$u(z)=P_{\alpha}[f_{\rho}](z).$$
 If $\alpha=0$, let  $z=re^{i\theta}$, and consider variable substitution of Section 2.5, then we have
\begin{align*}
\frac{\partial}{\partial z}u(z)&=\int^{2\pi}_{0}\frac{-\bar{z}|1-ze^{-it}|^{2}
+(1-r^{2})e^{-it}(1-\bar{z}e^{it})}{|1-ze^{-it}|^{4}}
f(e^{it})\frac{dt}{2\pi}\\
&=\frac{e^{-i\theta}}{2\pi}\int^{2\pi}_{0}\frac{-r\frac{(1-r^{2})^{2}}{|1+re^{-is}|^{2}}
+(1-r^{2})\left(\frac{r+e^{-is}}{1+re^{-is}}-r\right)}{\frac{(1-r^{2})^{4}}{|1+re^{-is}|^{4}}}
f(e^{it})\frac{1-r^{2}}{|1+re^{is}|^{2}}ds\\
&=\frac{1}{2\pi}(1-r^{2})^{-1}\int^{2\pi}_{0}e^{-i(\theta+s)}f(e^{it})ds\\
&=\frac{1}{2\pi}(1-r^{2})^{-1}\int^{2\pi}_{0}e^{-i(\theta+s)}(1-\rho^{2})^{\frac{2}{1-p}}(1+\cos s)^{\frac{1}{p-1}}e^{is}ds\\
&=\frac{e^{-i\theta}}{2\pi}(1-r^{2})^{-1}(1-\rho^{2})^{\frac{2}{1-p}}\int^{2\pi}_{0}(1+\cos s)^{\frac{1}{p-1}}ds.
\end{align*}
Let $r=\rho\rightarrow 1$, then
\begin{align*}
&\lim_{r=\rho\rightarrow 1}\frac{(1-r^{2})^{1+\frac{1}{p}}\,|u_{z}|}{\|f_{\rho}\|_{p}}
=2^{-1+\frac{2}{p}}\pi^{-\frac{1}{q}}
\left(\int^{2\pi}_{0}(1+\cos s)^{\frac{1}{p-1}}ds\right)^{\frac{1}{q}}\\
&=2^{1-\frac{2}{q}}\pi^{-\frac{1}{q}}
\left(\int^{2\pi}_{0}(1+\cos s)^{q-1}ds\right)^{\frac{1}{q}}
=\lim_{r\rightarrow 1}C_{0,p}(r)=C_{0,p}.
\end{align*}This means the constant $C_{0,p}$ is sharp.

The estimation of $|u_{\bar{z}}|$ is similar and the sharpness of constant  $C_{0,p}$  can be verified by taking $f_{\rho}(e^{it}):=(1-\rho^{2})^{\frac{2}{1-p}}(1+\cos s)^{\frac{1}{p-1}}e^{-is}$. We  omit the proof.

\end{proof}

\begin{proof}[\textbf{Proof of Theorem \ref{Th1.8}}]

(1) \quad Equation (\ref{3.17}) implies that

\begin{align}\label{3.40}
|u_{r}|\leq c_{\alpha}(1-r^{2})^{\alpha}\int^{2\pi}_{0}\frac{|2(\alpha+1)r|1-re^{i(\theta-t)}|^{2}+(\alpha+2)(1-r^{2})(r-\cos(\theta-t))|}
{|1-re^{i(\theta-t)}|^{4+\alpha}}|f(e^{it})|\frac{dt}{2\pi}.
\end{align}
Let \begin{align*}
I_{4}=\int^{2\pi}_{0}\frac{|2(\alpha+1)r|1-re^{i(\theta-t)}|^{2}+(\alpha+2)(1-r^{2})(r-\cos(\theta-t))|}
{|1-re^{i(\theta-t)}|^{4+\alpha}}\frac{dt}{2\pi}.
\end{align*}
Then by Jensen's inequality, we have that
\begin{align*}
|u_{r}|^{p}&\leq (c_{\alpha}(1-r^{2})^{\alpha})^{p}I_{4}^{p}\left(\int^{2\pi}_{0}\frac{|2(\alpha+1)r|1-re^{i(\theta-t)}|^{2}+(\alpha+2)(1-r^{2})(r-\cos(\theta-t))|}
{|1-re^{i(\theta-t)}|^{4+\alpha}I_{4}}|f(e^{it})|\frac{dt}{2\pi}\right)^{p}\nonumber\\
&\leq (c_{\alpha}(1-r^{2})^{\alpha})^{p}I_{4}^{p-1}\int^{2\pi}_{0}\frac{|2(\alpha+1)r|1-re^{i(\theta-t)}|^{2}+(\alpha+2)(1-r^{2})(r-\cos(\theta-t))|}
{|1-re^{i(\theta-t)}|^{4+\alpha}}|f(e^{it})|^{p}\frac{dt}{2\pi}.
\end{align*}
It follows that
\begin{align}\label{3.41}
&\quad\frac{1}{2\pi}\int^{2\pi}_{0}|u_{r}|^{p}d\theta\nonumber\\
&\leq (c_{\alpha}(1-r^{2})^{\alpha})^{p}I_{4}^{p-1}\int^{2\pi}_{0}
\left(\int^{2\pi}_{0}\frac{|2(\alpha+1)r|1-re^{i(\theta-t)}|^{2}+(\alpha+2)(1-r^{2})(r-\cos(\theta-t))|}
{|1-re^{i(\theta-t)}|^{4+\alpha}}|f(e^{it})|^{p}\frac{dt}{2\pi}\right)\frac{ d\theta}{2\pi}\nonumber\\
&=(c_{\alpha}(1-r^{2})^{\alpha})^{p}I_{4}^{p}\|f\|^{p}_{L^{p}(\mathbb{T})}
\end{align}
After charge of variables as in section 2.5,  direct computation leads to
 \begin{align*}
I_{4}=(1-r^{2})^{-1-\alpha}\int^{2\pi}_{0}\frac{|\alpha r-(\alpha+2)\cos s|}{|1+re^{-is}|^{-\alpha}}\frac{ds}{2\pi}.
\end{align*}Then (\ref{3.41}) implies that
\begin{align}\label{3.42}
M_{p}(r,u_{r})\leq \frac{\widetilde{D}_{\alpha}(r)}{1-r^{2}}\|f\|_{L^{p}(\mathbb{T)}},
\end{align}
where \begin{align*}
\widetilde{D}_{\alpha}(r)=c_{\alpha}\int^{2\pi}_{0}\frac{|\alpha r-(\alpha+2)\cos s|}{|1+re^{-is}|^{-\alpha}}\frac{ds}{2\pi}.
\end{align*}
Let \begin{equation}\label{3.43}
D_{\alpha}(r)=c_{\alpha}\left(\frac{ \alpha+2}{2\pi}\int^{2\pi}_{0}|\cos s||1+re^{-is}|^{\alpha}ds+\frac{|\alpha| r}{2\pi}\int^{2\pi}_{0}|1+re^{-is}|^{\alpha}ds\right)
\end{equation}and
\begin{equation} \label{3.44}\; D_{\alpha}= \left\{
\begin{array}{rr}
c_{\alpha}\left(\frac{ \alpha+2}{\pi}2^{\frac{\alpha}{2}-1}\int^{2\pi}_{0}|\cos s|(1+\cos s)^{\alpha/2}ds+|\alpha| \frac{\Gamma((\alpha+1)}{\left(\Gamma(\frac{\alpha+2}{2})\right)^{2}}\right), &\quad \alpha>2;\\
c_{\alpha}\left(\frac{ \alpha+2}{\pi}2^{\frac{\alpha}{2}-1}\int^{2\pi}_{0}|\sin s|(1+\cos s)^{\alpha/2}ds+|\alpha| \frac{\Gamma((\alpha+1)}{\left(\Gamma(\frac{\alpha+2}{2})\right)^{2}}\right),  &\quad -1<\alpha\leq 2.
\end{array}\right.
\end{equation}Considering Lemma \ref{lem2.3}, we have that $\widetilde{D}_{\alpha}(r)\leq D_{\alpha}(r)\leq D_{\alpha}$. Therefore,  (\ref{3.42}) implies (\ref{1.16}) holds with $D_{\alpha}(r)\leq D_{\alpha}$.

Specially, if $\alpha=0$, then $D_{0}(r)=D_{0}=\frac{4}{\pi}$.

Next, let us show the constant $D_{0}$  is sharp.  Let $0<\rho<1$, and define
\begin{align}
f_{\rho}(e^{it})=sign(\cos s),
\end{align}where $s$ and $t$ satisfy the relation (\ref{3.25}).
 Considering the variable change of (\ref{3.25}), we have \begin{align*}&\quad \frac{1}{2\pi}\int^{2\pi}_{0}|f_{\rho}(e^{it})|^{p}dt=\frac{1}{2\pi}\int^{2\pi}_{0}\frac{1-\rho^{2}}{|1+\rho e^{is}|^{2}}ds=1.
\end{align*}
Thus, $\|f_{\rho}\|_{L^{p}(\mathbb{T})}=1$.

 Take $$u(z)=P_{\alpha}[f_{\rho}](z)$$with $\alpha=0$.
Let $r=\rho$, we have
\begin{align*}
\frac{\partial}{\partial r}u(re^{i\theta})&=-\frac{1}{2\pi}\int^{2\pi}_{0}\frac{2r|1-r\frac{r+e^{-is}}{1+r e^{-is}}|^{2}-2(1-r^{2})\frac{(1-r^{2})(r+\cos s)}{|1+re^{-is}|^{2}}}
{|1-r\frac{r+e^{-is}}{1+r e^{-is}}|^{4}}f(e^{it})\frac{1-r^{2}}{|1+r e^{-is}|^{2}}ds\\
&=\frac{1}{\pi}\int^{2\pi}_{0}\frac{\cos s}{1-r^{2}}f(e^{it})ds
=\frac{1}{\pi}\int^{2\pi}_{0}\frac{|\cos s|}{1-r^{2}}ds
=\frac{4}{\pi(1-r^{2})}.
\end{align*}
It follows that $M_{p}(r,u_{r})=\frac{4}{\pi(1-r^{2})}=\frac{D_{0}}{(1-r^{2})}\|f_{\rho}\|_{p}$.
This shows  the constant $D_{0}$ is sharp.

(2)\quad Equation (\ref{3.26}) implies that
 \begin{align}\label{3.46}
|u_{\theta}|\leq\frac{(\alpha+2)c_{\alpha}}{2\pi}(1-r^{2})^{\alpha+1}\int^{2\pi}_{0}\frac{r|\sin(\theta-t)|}
{(1+r^{2}-2r\cos(\theta-t))^{2+\alpha/2}}|f(e^{it})|dt
\end{align}

Let  \begin{align*}
I_{5}=\int^{2\pi}_{0}\frac{r|\sin(\theta-t)|}
{(1+r^{2}-2r\cos(\theta-t))^{2+\alpha/2}}dt.
\end{align*}Then direct computation leads to that
  \begin{align}\label{3.47}
I_{5}&=2\int^{\pi}_{0}\frac{r\sin(\theta-t)}
{(1+r^{2}-2r\cos(\theta-t))^{2+\alpha/2}}dt\nonumber\\
&=\int^{\pi}_{0}(1+r^{2}-2r\cos T))^{-{(2+\alpha/2)}}d(1+r^{2}-2r\cos T)\nonumber\\
&=\frac{2}{\alpha+2}\left(1-\left(\frac{1-r}{1+r}\right)^{\alpha+2}\right)\frac{1}{(1-r)^{\alpha+2}}.
\end{align}
It follows from (\ref{3.46}) and Jensen's inequality that
 \begin{align*}
|u_{\theta}|^{p}&\leq\left(\frac{(\alpha+2)c_{\alpha}}{2\pi}(1-r^{2})^{\alpha+1}\right)^{p}
\left(I_{5}\cdot\int^{2\pi}_{0}\frac{r|\sin(\theta-t)|}
{(1+r^{2}-2r\cos(\theta-t))^{2+\alpha/2}I_{5}}|f(e^{it})|dt\right)^{p}\nonumber\\
&\leq\left(\frac{(\alpha+2)c_{\alpha}}{2\pi}(1-r^{2})^{\alpha+1}\right)^{p}I_{5}^{p-1}\cdot\int^{2\pi}_{0}\frac{r|\sin(\theta-t)|}
{(1+r^{2}-2r\cos(\theta-t))^{2+\alpha/2}}|f(e^{it})|^{p}dt\nonumber.
\end{align*}
Therefore,
\begin{align}\label{3.48}
&\quad\frac{1}{2\pi}\int^{2\pi}_{0}|u_{\theta}|^{p} d\theta\nonumber\\
&\leq\left(\frac{(\alpha+2)c_{\alpha}}{2\pi}(1-r^{2})^{\alpha+1}\right)^{p}I_{5}^{p-1}\cdot\frac{1}{2\pi}\int^{2\pi}_{0}\left(\int^{2\pi}_{0}\frac{r|\sin(\theta-t)|}
{(1+r^{2}-2r\cos(\theta-t))^{2+\alpha/2}}|f(e^{it})|^{p}dt\right)d\theta\nonumber\\
&\leq\left(\frac{(\alpha+2)c_{\alpha}}{2\pi}(1-r^{2})^{\alpha+1}\right)^{p}I_{5}^{p}\|f\|^{p}_{L^{p}(\mathbb{T})}.
\end{align}
Let\begin{align}\label{3.49}
E_{\alpha}(r)=\frac{c_{\alpha}}{\pi}((1+r)^{\alpha+2}-(1-r)^{\alpha+2}).
\end{align}Then it follows from (\ref{3.47}) and (\ref{3.48})  that
\begin{align*}
M_{p}(r, u_{\theta})\leq\frac{E_{\alpha}(r)}{1-r^{2}}\|f\|_{L^{p}(\mathbb{T})}.
\end{align*}
Obviously, \begin{align}\label{3.50}
E_{\alpha}(r)\leq  c_{\alpha}2^{\alpha+2}/\pi=:E_{\alpha}
\end{align} for $r\in [0,1)$.

Next, let us show the constant $E_{0}=\frac{4}{\pi}$  is sharp.  Let $0<\rho<1$, and define
\begin{align}\label{3.51}
f_{\rho}(e^{it})=sign(\sin s),
\end{align}where $s$ and $t$ satisfy the relation (\ref{3.25}).
Then  we have $\|f_{\rho}\|_{L^{p}(\mathbb{T})}=1$.

 Take $$u(z)=P_{\alpha}[f_{\rho}](z).$$
If $\alpha=0$,
taking $r=\rho$, (\ref{3.26}) reduces to
\begin{align*}
\frac{\partial}{\partial \theta}u(re^{i\theta})&=-\frac{1}{\pi}r(1-r^{2})\int^{2\pi}_{0}\frac{\sin(\theta-t)}
{|1-re^{i(\theta-t)}|^{4}}f(e^{it})dt
=\frac{1}{\pi}r(1-r^{2})^{-1}\int^{2\pi}_{0}\sin s \,f(e^{it})ds
=\frac{4r}{\pi(1-r^{2})}.
\end{align*}
It follows that \begin{align*}\lim_{\rho=r\rightarrow 1}\frac{(1-r^{2})M_{p}(r,u_{\theta})}{\|f_{\rho}\|_{L^{p}(\mathbb{T})}}=\frac{4}{\pi}=E_{0}.
\end{align*}
This shows that the constant $E_{0}$ is sharp.

(3) It follows from (\ref{3.33}) that
\begin{align}\label{3.52}
|u_{z}|&\leq c_{\alpha}(1-r^{2})^{\alpha}\int^{2\pi}_{0}\frac{|-(\alpha+1)\bar{z}|1-ze^{-it}|^{2}
+(1+\alpha/2)(1-r^{2})e^{-it}(1-\bar{z}e^{it})|}{|1-ze^{-it}|^{\alpha+4}}
|f(e^{it})|\frac{dt}{2\pi}.
\end{align}
Let \begin{align*}
I_{6}=\int^{2\pi}_{0}\frac{|-(\alpha+1)\bar{z}|1-ze^{-it}|^{2}
+(1+\alpha/2)(1-r^{2})e^{-it}(1-\bar{z}e^{it})|}{|1-ze^{-it}|^{\alpha+4}}\frac{dt}{2\pi}.
\end{align*}
After change of variables as in section 2.5, we have
 \begin{align*}
 I_{6}=(1-r^{2})^{-1-\alpha}\int^{2\pi}_{0}\frac{|-\alpha r/2
+(1+\alpha/2)e^{-is}|}{|1+re^{-is}|^{-\alpha}}\frac{ds}{2\pi}.
\end{align*}
Using Jensen's inequality, from  (\ref{3.52}), we get
\begin{align*}
|u_{z}|^{p}&\leq (c_{\alpha}(1-r^{2})^{\alpha})^{p}\left(I_{6}\cdot\int^{2\pi}_{0}\frac{|-(\alpha+1)\bar{z}|1-ze^{-it}|^{2}
+(1+\alpha/2)(1-r^{2})e^{-it}(1-\bar{z}e^{it})|}{|1-ze^{-it}|^{\alpha+4}I_{6}}
|f(e^{it})|\frac{dt}{2\pi}\right)^{p}\nonumber\\
&\leq (c_{\alpha}(1-r^{2})^{\alpha})^{p}I_{6}^{p-1}\int^{2\pi}_{0}\frac{|-(\alpha+1)\bar{z}|1-ze^{-it}|^{2}
+(1+\alpha/2)(1-r^{2})e^{-it}(1-\bar{z}e^{it})|}{|1-ze^{-it}|^{\alpha+4}}
|f(e^{it})|^{p}\frac{dt}{2\pi}.
\end{align*}
It follows that
\begin{align}\label{3.53}
&\quad\frac{1}{2\pi}\int^{2\pi}_{0}|u_{z}|^{p}d\theta\nonumber\\
 &\leq (c_{\alpha}(1-r^{2})^{\alpha})^{p}I_{6}^{p-1}\int^{2\pi}_{0}
 \left(\int^{2\pi}_{0}\frac{|-(\alpha+1)\bar{z}|1-ze^{-it}|^{2}
+(1+\alpha/2)(1-r^{2})e^{-it}(1-\bar{z}e^{it})|}{|1-ze^{-it}|^{\alpha+4}}
|f(e^{it})|^{p}\frac{dt}{2\pi}\right)\frac{ d\theta}{2\pi}\nonumber\\
&\leq (c_{\alpha}(1-r^{2})^{\alpha})^{p}I_{6}^{p}\|f\|^{p}_{L^{p}(\mathbb{T})}.
\end{align}
Let \begin{align*}
\widetilde{F}_{\alpha}(r)=c_{\alpha}\int^{2\pi}_{0}\frac{|-\alpha r/2
+(1+\alpha/2)e^{-is}|}{|1+re^{-is}|^{-\alpha}}\frac{ds}{2\pi}.
\end{align*}Then (\ref{3.53}) reduces to
\begin{align*}
M_{p}(r, u_{z})\leq \frac{\widetilde{F}_{\alpha}(r)}{1-r^{2}}\|f\|^{p}_{L^{p}(\mathbb{T})}.
\end{align*}
 Let
\begin{equation}\label{3.54}
F_{\alpha}(r)=\frac{c_{\alpha}(|\alpha |r/2
+1+\alpha/2)}{2\pi}\int^{2\pi}_{0}\frac{ds}{|1+re^{-is}|^{-\alpha}}.
\end{equation}
and \begin{equation}\label{3.55}
F_{\alpha}= c_{\alpha}(|\alpha |/2
+1+\alpha/2)\frac{\Gamma((\alpha+1)}{\left(\Gamma(\frac{\alpha+2}{2})\right)^{2}}.
\end{equation} It follows that $\widetilde{F}_{\alpha}(r)\leq F_{\alpha}(r)\leq F_{\alpha}$.  Therefore,  (\ref{1.18}) holds with $F_{\alpha}(r)\leq F_{\alpha}$.

If $\alpha=0$, then $F_{0}(r)=F_{0}=1$. Next, we show the constant $F_{0}$ is sharp.

 Let $0<\rho<1$, and define
\begin{align}
f_{\rho}(e^{it})=e^{is},
\end{align}where $s$ and $t$ satisfy the relation of (\ref{3.25}).
Then we have
$\|f_{\rho}\|_{L^{p}(\mathbb{T})}=1$.

If $\alpha=0$, let  $z=re^{i\theta}$, then we have
\begin{align*}
\frac{\partial}{\partial z}u(z)&=\int^{2\pi}_{0}\frac{-\bar{z}|1-ze^{-it}|^{2}
+(1-r^{2})e^{-it}(1-\bar{z}e^{it})}{|1-ze^{-it}|^{4}}
f(e^{it})\frac{dt}{2\pi}\\
&=\frac{1}{2\pi}(1-r^{2})^{-1}\int^{2\pi}_{0}e^{-i(\theta+s)}f(e^{it})ds=(1-r^{2})^{-1}e^{-i\theta}.
\end{align*}
Let $r=\rho\rightarrow 1$, then
 \begin{align*}
&\lim_{r=\rho\rightarrow 1}\frac{(1-r^{2})\,M_{p}(r, u_{z})}{\|f_{\rho}\|_{L^{p}(\mathbb{T})}}=F_{0},
\end{align*}which means the constant $F_{0}$ is sharp.

The estimation of $M_{p}(r,u_{\bar{z}})$ is similar. The sharpness of constant  $F_{0}$  can be verified by taking  $f_{\rho}(e^{it})=e^{-is}$. We  omit the proof.
\end{proof}

\begin{proof}[\textbf{Proof of Theorem \ref{Th1.9}}]
Suppose $f(e^{it})$ is the boundary function of $\alpha$-harmonic function $u(z)$. Then
\begin{equation*}u(z)=P_{\alpha}[f](z)=\frac{1}{2\pi}\int^{2\pi}_{0}K_{\alpha}(r,\theta-t)f(e^{it})dt,\quad z=re^{i\theta}\in\mathbb{D}.
\end{equation*}It follows that
\begin{align}\label{3.57}
\frac{\partial}{\partial\theta}u(re^{i\theta}))
&=\frac{1}{2\pi}\int^{2\pi}_{0}\frac{\partial}{\partial\theta}\left(K_{\alpha}(r,\theta-t)\right)f(e^{it})dt\nonumber\\
&=-\frac{1}{2\pi}\int^{2\pi}_{0}f(e^{it})\frac{\partial}{\partial t}K_{\alpha}(r,\theta-t)dt\nonumber\\
 &=\frac{1}{2\pi}\int^{2\pi}_{0}K_{\alpha}(r,\theta-t)df(e^{it}),
\end{align}
because the function $t\mapsto f(e^{it})$ is of bounded variation on $[0,2\pi]$; the last integral is regarded as the Stieltjes one. Then by (\ref{3.57}),
\begin{align}\label{3.58}
 \frac{\partial}{\partial\theta}u(re^{i\theta}))
=\frac{1}{2\pi}\lim_{n\rightarrow\infty}\sum^{n}_{k=1}K_{\alpha}(r,\theta-2k\pi/n)(f(e^{2k\pi i/n})-f(e^{i2(k-1)\pi i/n})),
\end{align}
Hence, applying Fatou's limiting integral lemma, we obtain
\begin{align}
\int^{2\pi}_{0}\left|\frac{\partial}{\partial\theta}u(re^{i\theta}))\right|d\theta
&=\frac{1}{2\pi}\int^{2\pi}_{0}\left|\lim_{n\rightarrow\infty}\sum^{n}_{k=1}P_{p,q}(r,\theta-2k\pi/n)(f(e^{2k\pi i/n})-f(e^{i2(k-1)\pi i/n}))\right|d\theta\nonumber\\
&\label{3.59}\leq\frac{1}{2\pi}\liminf_{n\rightarrow \infty}\int^{2\pi}_{0}\left|\sum^{n}_{k=1}P_{p,q}(r,\theta-2k\pi/n)(f(e^{2k\pi i/n})-f(e^{i2(k-1)\pi i/n}))\right|d\theta.
\end{align}
Since
\begin{align*}
&\frac{1}{2\pi}\int^{2\pi}_{0}\left|\sum^{n}_{k=1}P_{p,q}(r,\theta-2k\pi/n)(f(e^{2k\pi i/n})-f(e^{i2(k-1)\pi i/n}))\right|d\theta\nonumber\\
&\leq \sum^{n}_{k=1}\left[\left|f(e^{2k\pi i/n})-f(e^{2(k-1)\pi i/n})\right|\frac{1}{2\pi}\int^{2\pi}_{0}|K_{\alpha}(r,\theta-2k\pi/n)|d\theta\right]\nonumber\\
&\leq c_{\alpha}F(-\alpha/2, -\alpha/2: 1:r^{2})\sum^{n}_{k=1}\left|f(e^{2k\pi i/n})-f(e^{2(k-1)\pi i/n})\right|\nonumber\\
&\leq  c_{\alpha}F(-\alpha/2, -\alpha/2: 1:r^{2})|\gamma|\nonumber
\end{align*}and for $z=re^{i\theta}$,
\begin{align}\label{3.60}
\frac{\partial}{\partial\theta}u(re^{i\theta}))=i[z\partial u(z)-\bar{z}\bar{\partial}u(z)],
\end{align}we conclude from (\ref{3.59}) that
\begin{align}\label{3.61}
\int_{\mathbb{T}_{r}}(|\partial u(z)|-|\bar{\partial}u(z)|)|dz|
\leq\int^{2\pi}_{0}(|z\partial u(z)-\bar{z}\bar{\partial}u(z)|)d\theta\leq  c_{\alpha}F(-\alpha/2, -\alpha/2: 1:r^{2})|\gamma|.
\end{align}
By the assumption, the mapping $u$ is $K$-quasiconformal, which means that
\begin{align}\label{3.62}
(K+1)|\bar{\partial}u(z)|\leq(K-1)|\partial u(z)|, \quad z\in\mathbb{D}.
\end{align}Hence by (\ref{3.61}),
\begin{align}\label{3.63}
\int_{\mathbb{T}_{r}}(|\bar{\partial}u(z)|+|\partial u(z)|)|dz|
\leq K\int_{\mathbb{T}_{r}}(|\partial u(z)|-|\bar{\partial}u(z)|)|dz|
\leq K  c_{\alpha}F(-\alpha/2, -\alpha/2: 1:r^{2})|\gamma|.
\end{align}
Combining this with (\ref{3.61}) and (\ref{3.63}) leads to
\begin{align}\label{3.64}
\int_{\mathbb{T}_{r}}|\partial u(z)||dz|
\leq \frac{K+1}{2}  c_{\alpha}F(-\alpha/2, -\alpha/2: 1:r^{2})|\gamma|.
\end{align}and
\begin{align}\label{3.65}
\int_{\mathbb{T}_{r}}|\bar{\partial} u(z)||dz|
\leq \frac{K-1}{2} c_{\alpha}F(-\alpha/2, -\alpha/2: 1:r^{2})|\gamma|.
\end{align}
By Lemma \ref{lem2.2}, $F(-\alpha/2,-\alpha/2;1;r^{2})$ is increasing on $r\in (0,1)$. Then by (\ref{3.64}) and (\ref{3.65}) , we can get
\begin{align*}
\int_{\mathbb{T}_{r}}|\partial u(z)||dz|
\leq \frac{K+1}{2} |\gamma|.
\end{align*}and
\begin{align*}\int_{\mathbb{T}_{r}}|\bar{\partial} u(z)||dz|
\leq \frac{K-1}{2} |\gamma|.
\end{align*}Therefore, it leads to (\ref{1.21}) and (\ref{1.22}).

Obviously, when $K\rightarrow1$ and $\alpha=0$, we take $u(z)=z$, then the equality holds.
\end{proof}

 \begin{proof}[\textbf{Proof of Theorem \ref{Th1.10}}]
 It follows (\ref{1.22}) that $\Omega$ is bounded by a rectifiable Jordan cure $\gamma$. Hence, the function $t\rightarrow f(e^{it})$ is of bounded variation on $[0, 2\pi]$ and as in the proof of Theorem \ref{Th1.9}, the equality (\ref{3.58}) holds.  It follows from (\ref{1.23}) that for all $n\in \mathbb{N}$ and $k=1, 2,...n$,
 \begin{equation}\label{3.66}
|f(e^{2\pi k i/n})-f(e^{2\pi (k-1)i/n})|\leq L |e^{2\pi k i/n}-e^{2\pi (k-1)i/n}|=2L \sin\frac{\pi}{n}.
\end{equation}
 We conclude from (\ref{3.58}) and  (\ref{3.66}) that for every $\theta\in \mathbb{R}$
  \begin{align}\label{3.67}
\left|\frac{\partial}{\partial \theta}u(re^{i\theta})\right|
&\leq \lim_{n\rightarrow \infty}\left[\frac{1}{2\pi}\sum^{n}_{k=1}|K_{\alpha}(ze^{-i\frac{2k\pi}{n}})|\cdot2L \sin\frac{\pi}{n}\right|\nonumber\\
&=L\lim_{n\rightarrow\infty}\frac{1}{2\pi}\sum^{n}_{k=1}\frac{2\pi}{n}|K_{\alpha}(ze^{-i\frac{2k\pi}{n}})|\nonumber\\
&=L\cdot\frac{1}{2\pi}\int^{2\pi}_{0}|K_{\alpha}(ze^{-it})|dt\nonumber\\
&=Lc_{\alpha}F(-\alpha/2,-\alpha/2;1; r^{2}).
 \end{align}
Since for $z=re^{i\theta}$ the equality (\ref{3.60}) holds, we conclude from (\ref{3.67}) that
 \begin{align}\label{3.68}
 r(|\partial u|-|\bar{\partial}u|)\leq |z\partial u(z)-\bar{z}\bar{\partial}u(z)|\leq Lc_{\alpha}F(-\alpha/2,-\alpha/2;1; r^{2}).
\end{align}
By the assumption, the mapping $u(z)$ is $K$-qc., which means that (\ref{3.62}) holds. Here by (\ref{3.68}),
 \begin{align}\label{3.69}
 r(|\partial u|+|\bar{\partial}u|)\leq rK(|\partial u|-|\bar{\partial}u|)\leq K Lc_{\alpha}F(-\alpha/2,-\alpha/2;1; r^{2}).
\end{align}
Combining the inequalities (\ref{3.68}) and (\ref{3.69}), we obtain
\begin{align*}
|z\partial u|\leq \frac{K+1}{2}Lc_{\alpha}F(-\alpha/2,-\alpha/2;1; r^{2})\leq\frac{K+1}{2}L
\end{align*}and
\begin{align*}
|\bar{z}\bar{\partial}u|\leq \frac{K-1}{2}Lc_{\alpha}F(-\alpha/2,-\alpha/2;1; r^{2})\leq\frac{K-1}{2}L.
\end{align*}

Obviously, when $K\rightarrow1$ and $\alpha=0$, we take $u(z)=z$, then the equality holds.
 \end{proof}

\medskip

{\bf Declarations}\\

{\bf Conflict of interests}\quad The authors declare that they have no conflict of interest.\\

{\bf Data availability statement}\quad This manuscript has no associated date.

\medskip

\end{document}